 0  0

\documentclass[a4paper]{amsart}
\usepackage{amssymb,graphics,wrapfig,caption,ifthen,array}
\usepackage[all,knot]{xy}
\xyoption{import,arc}

\newcommand{\ad} [1] { {\rm ad}_{#1} }
\newcommand{\Ad} [1] { {\rm Ad}_{#1} }

\newcommand{\R} [0] { {\bf{R}} }
\newcommand{\T} [0] { {\bf{T}} }
\newcommand{\Z} [0] { {\bf{Z}} }

\newcommand{\Q} [0] { {\bf{Q}} }
\newcommand{\C} [0] { {\bf{C}} }

\newcommand \im { {\rm im}\, }

\newcommand {\h} [0] { {\mathfrak h} }
\newcommand {\z} [0] { {\mathfrak z} }
\newcommand {\g} [0] { {\mathfrak g} }
\newcommand {\A} [0] { {\mathfrak a} }
\def\B{{\mathfrak b}}

\def\mc{{\mathfrak m}}
\def\d{{\mathrm d}}
\def\nw#1{\Omega(#1)}
\def\forms#1#2{\wedge^{#1}(#2)}
\def\U{{\rm U}}
\def\SU{{\rm SU}}

\def\u{{\rm u}}
\def\su{{\rm su}}
\def\rank{{\rm rank\, }}
\def\diag{{\rm diag}}
\def\mm{{\Psi}}

\def\L{{\mathfrak L}}

\def\Cw#1{{\rm C}^{\omega}{\rm (}#1{\rm )}}

\def\uu{{\mathfrak u}}
\def\ddim#1{{\rm ddim(}#1{\rm )}}
\def\drank#1{{\rm drank(}#1{\rm )}}
\def\Ch#1{ C^{\infty}(#1) }
\def\pb#1#2{{ \left\{ #1, #2  \right\} }}
\def\spn#1{{\rm  span} \left\{ #1  \right\}}
\def\orbit#1{{\mathcal O}_{#1}}
\def\t{{\mathfrak t}}
\def\ie{{\it i.e.\ }}
\def\cf{{\it c.f.\ }}
\def\jj{{\textrm{\j}}}

\def\Deck{{\Delta}}

\def\Trace#1{\textrm{Trace}\, #1}

\def\x{{\bf x}}
\def\y{{\bf y}}
\def\w{{\bf w}}
\def\z{{\bf z}}
\def\re#1{{\rm Re\,} #1}
\def\g{{\bf g}}
\def\bg{{\bf \bar{g}}}
\def\L{{\mathsf G}}
\def\adj#1{\textrm{adj}(#1)}
\def\S{{\mathsf S}}
\def\V{{\mathsf V}}
\def\ul{{\,\underline{\ }\hspace{-1.5mm}\lambda}}
\def\bl{{\bar{\lambda}}}
\newcommand{\filename}[2]{\ifthenelse{\equal{#1}{1}}{$ $\newline {\tt #2} }{}}
\def\H{{\mathsf H}}
\def\O{{\rm O}}
\def\SO{{\rm SO}}

\newtheorem{thm}{Theorem}
\newtheorem{qthm}{Theorem}
\newtheorem{lemma}{Lemma}[section]
\newtheorem{prop}{Proposition}[section]
\newtheorem{definition}{Definition}
\newtheorem{corollary}{Corollary}
\newtheorem{question}{Question}

\newcounter{remark}
\renewcommand{\theremark}{\Roman{remark}}
\newenvironment{remark}{\begin{trivlist}\item[]\refstepcounter{remark}%
	{\tt Remark\ \theremark.}\nobreak\noindent\rm\ignorespaces}{%
	\end{trivlist}}
\newcounter{claim}
\renewcommand{\theclaim}{\arabic{thm}.\arabic{claim}}
\newenvironment{claim}{
\refstepcounter{claim}%
\medskip
\noindent{\bf Claim\ \theclaim .}%
\hspace{2mm}\nobreak\noindent\it\ignorespaces}{}
\newcounter{chec}

\newenvironment{chec}{
\refstepcounter{chec}%
\medskip
\noindent
{\bf Check.}%
\hspace{2mm}\nobreak\noindent\rm\ignorespaces
}
{\hfill$\diamondsuit$}

\begin{document}

\title[The Maslov Cocycle and Integrability]{The Maslov cocycle,
  smooth structures and real-analytic complete integrability
\filename{}{maslov.tex $->$ pi1II.tex}}

\author{Leo T. Butler}
\address{
School of Mathematics, The University of Edinburgh, Edinburgh, UK, EH9
3JZ }
\email{l.butler@ed.ac.uk}
%

\begin{abstract} 
This paper proves two main results. First, it is shown that if
$\Sigma$ is a smooth manifold homeomorphic to the standard $n$-torus
$\T^n = \R^n/\Z^n$ and $H$ is a real-analytically completely
integrable convex hamiltonian on $T^* \Sigma$, then $\Sigma$ is
diffeomorphic to $\T^n$. Second, it is proven that for some
topological $7$-manifolds, the cotangent bundle of each smooth
structure admits a real-analytically completely integrable riemannian
metric hamiltonian. A version of this paper appears in {\em American
Journal of Mathematics}. 131(5) : 1311--1336, October 2009.
\end{abstract}
\date{\today}
\maketitle

\section{Introduction}
One of the most intriguing facts of differential topology is that a
topological manifold may admit several distinct smooth structures. An
important smooth invariant of a smooth manifold is the cotangent
bundle, so a smooth dynamical system on the cotangent bundle ought, in
principle, to reflect the smooth structure. In the present note, it is
shown that the existence of a real-analytically integrable convex
hamiltonian on the cotangent bundle is a non-trivial smooth invariant.

\subsection{Complete integrability} \label{ssec:ci} 
The cotangent bundle of a smooth manifold $\Sigma$ admits a canonical
symplectic form $\omega=\sum \d y_i \wedge \d x^i$, where $x^i$ are
coordinates on $\Sigma$ and $y_i$ are the induced fibre
coordinates. A symplectic form permits one to define a Poisson algebra
structure on $C^{\infty}(T^*\Sigma)$ and consequently each smooth
function $H : T^* \Sigma \to \R$ induces a hamiltonian vector field
$X_H$ defined by
\begin{equation} \label{eq:xh}
X_H = \pb{H}{\ } \qquad \implies X_H= \left\{ \begin{array}{lcl}
\dot{x}^i &=& \frac{\partial H}{\partial y_i},\\
\dot{y}_i &=& -\frac{\partial H}{\partial x^i}.
\end{array}  
\right.
\end{equation}
A first integral of the hamiltonian vector field $X_H$ is a smooth
function $F$ which Poisson commutes with $H$: $\pb{H}{F} = 0$. If
$X_H$ has $n=\dim \Sigma$ functionally independent first integrals
$F_1,\ldots,F_n$, and the first integrals pairwise Poisson commute,
then the compact regular level sets $\left\{ F_1=c_1,\ldots,\right.$
$F_n\left.=c_n \right\}$ are $n$-dimensional lagrangian tori and the
flow of $X_H$ is translation-type. In this case, one says that $X_H$ 
is {\em completely integrable}; if the first integrals are
real-analytic, one says that $X_H$ is real-analytically completely
integrable.

\subsection{Geometric semisimplicity} 
Let us abstract the notion of complete integrability. A smooth flow
$\varphi : M \times \R \to M$ is {\em integrable} if there is an open,
dense subset $R \subset M$ that is covered by angle-action charts
which conjugate $\varphi$ to a translation-type flow on the tori of
$\T^k \times \R^l$. There is an open dense subset $L \subset R$ fibred
by $\varphi$-invariant tori; let $f : L \to B$ be the induced smooth
quotient map and let $\Gamma = M - L$ be the {\em singular set}. If
$\Gamma$ is a tamely-embedded polyhedron, then $\varphi$ is said to be
$k$-{\em semisimple} with respect to $(f,L,B)$, or just
semisimple~\cite{Butler:2005}. Of most interest is when $\varphi$ is a
hamiltonian flow on a cotangent bundle or possibly a regular
iso-energy surface.
\begin{definition}[\cf \cite{Taimanov,Butler:2005}] \label{def:gss}
A hamiltonian flow is {\em geometrically semisimple} if it is
semisimple with respect to $(f,L,B)$ and $f$ is a lagrangian
fibration. It is {\em finitely geometrically semisimple} if, in
addition, each component of $B$ has a finite fundamental group.
\end{definition}
In this case, the lagrangian-ness of the fibres of $f$ implies that
$\varphi$ is completely integrable, so geometric semisimplicity is a
topologically-tame type of complete
integrability. Taimanov~\cite{Taimanov} introduced a related notion of
geometric simplicity, see sections 2.2-2.3 of~\cite{Butler:2005} for
further discussion. If $\varphi$ is real-analytically completely
integrable, then the triangulability of real-analytic sets implies
that $\varphi$ is finitely geometrically semisimple; in fact, in this
case $B$ may be taken to be a disjoint union of open balls. On the
other hand, geometric semisimplicity is a weaker property than
real-analytic complete integrability~\cite{Butler:2005}. A basic
question is:

\begin{question} \label{q:a}
What are the obstructions to the existence of a geometrically
semisimple (resp. semisimple, completely integrable) flow?
\end{question}

\subsection{Main Results} 
Recall that a topological $n$-torus is a topological manifold that is
homeomorphic to the standard $n$-torus $\T^n=\R^n/\Z^n$. An exotic
$n$-torus is a topological $n$-torus that is not diffeomorphic to
$\T^n$. Exotic $n$-tori may be constructed by connect summing with
exotic spheres, but not all arise this way.

\begin{thm} \label{thm:X}
If $\Sigma$ is an exotic $n$-torus, then there are no {\em finitely
geometrically semisimple} convex hamiltonians on $T^* \Sigma$. In
particular, there are no real-analytically completely integrable
convex hamiltonians on $T^* \Sigma$.
\end{thm}

The obstruction here is the {\em smooth} structure of the
configuration space. This is the first result that shows that a smooth
invariant may preclude real-analytic complete integrability; as such,
it prompts several questions.

\begin{question} \label{q:b}
If $\Sigma$ is an exotic torus, does
  there exist a completely integrable convex hamiltonian on
  $T^*\Sigma$?
\end{question}

Theorem \ref{thm:nope} in section \ref{sec:pe} strengthens Theorem
\ref{thm:X} by showing that there are no completely integrable
riemannian metrics on an exotic torus that are completely integrable
via a geodesic equivalence, see also \cite[Theorem 8]{Matveev}. This
suggests that the answer to Question
\ref{q:b} may be {\em no}.

It is important to note that the definitiveness of Theorem~\ref{thm:X}
is not general and may be atypical. The Gromoll-Meyer exotic
$7$-sphere is a biquotient Sp(2)//Sp(1) and so it inherits a
submersion metric from the bi-invariant metric on Sp(2)
\cite{GM}. Paternain and Spatzier \cite{PS1} proved the real-analytic
complete integrability of the geodesic flow of this submersion
metric. On the other hand, the remaining $12$ unoriented
diffeomorphism classes of the $7$-sphere are not known to possess such
geodesic flows.

\begin{question} \label{q:c}
Do all exotic $7$-spheres admit a real-analytically completely
  integrable convex hamiltonian?
\end{question}

And, more generally,

\begin{question} \label{q:d}
What are the smooth obstructions to the existence of a geometrically
  semisimple convex hamiltonian?
\end{question}

While the present paper does not answer Questions
\ref{q:b}--\ref{q:d}, it is able to answer the question for some
classes of topological $7$-manifolds with more than one smooth
structure. A {\em Witten-Kreck-Stolz space} $M_{k,l}$ is the smooth
$7$-manifold obtained by quotienting $S^5 \times S^3$ by the action of
$\U_1$ given by the representation $z \mapsto z^k \cdot I \oplus z^l
\cdot I : \U_1 \to \U_3 \oplus \U_2$, where $k$ and $l$ are coprime
integers. Kreck and Stolz showed that $M_{k,l}$ has a maximum of $28$
smooth structures; and, with modest conditions on $k$ and $l$, this
maximum is attained and each smooth structure is represented by some
$M_{k',l'}$~\cite{KS2}. This paper uses the work of Mykytyuk and
Panasyuk~\cite{MP} to show that

\begin{thm} \label{thm:KS}
There is a real-analytically completely integrable convex hamiltonian
on the cotangent bundle of each Witten-Kreck-Stolz space. In
particular, if $l=0 \bmod 4$, $l=0,3,4 \bmod 7$, $l \neq 0$ and
$\gcd(k,l)=1$, then each one of the $28$ diffeomorphism classes of
$M_{k,l}$ is the configuration space of a real-analytically completely
integrable convex hamiltonian.
\end{thm} 

The convex hamiltonian in all cases may be taken to be the hamiltonian
induced by the round metrics on $S^5$ and $S^3$. For each
Witten-Kreck-Stolz space, there is an $S^1$ fibre bundle
$S^1 \hookrightarrow M_{k,l} \to \C P^2 \times \C P^1$. In \cite{BJ2}
(resp. \cite{BJ1}) Bolsinov and Jovanovi\'c prove, {\it inter alia},
the real-analytic non-commutative (resp. complete) integrability of
the geodesic flows of certain homogeneous metrics on $\C P^2 \times \C
P^1$, see especially \cite[Remark 3.4]{BJ1}.

The present paper finishes by proving a similar result for the
Eschenburg and Aloff-Wallach $7$-manifolds. These manifolds are
obtained through a quotient of $\SU_3$ by a subgroup $V \cong \U_1$ of
the maximal torus of $\SU_3 \times \SU_3$. The existence of
real-analytically completely integrable geodesic flows on some special
Eschenburg spaces was proven by Paternain \& Spatzier and
Bazaikin~\cite{PS1,Baz}. The results of the present paper extend their
work. Kruggel \cite{Kr} has obtained a complete list of invariants
that classify the smooth structures on most Eschenburg spaces. It is
unknown if each topological Eschenburg space admits the maximum 28
smooth structures and each smooth structure is itself an Eschenburg
space. Numerical computations \cite{CEZ,Butler:2008a} suggest this may
be true for some families of Eschenburg spaces.

\subsubsection{Related work} \label{sssec:otherwork} 
Bialy and Polterovich \cite[Theorem 1.1]{BP2} prove that if $F \subset
T^*\T^2$ is an essential lagrangian torus that is invariant under a
convex hamiltonian flow, and without periodic points, then the natural
map $F \to \T^2$ is a diffeomorphism. Theorem \ref{thm:X} is based on
a generalization of their theorem to higher dimensions, see
Proposition \ref{prop:x3} and Remark \ref{rem:c} below.

Taimanov \cite{Taimanov} has proven that if a compact manifold
$\Sigma$ admits a real-analytically completely integrable geodesic
flow, then $\pi_1(\Sigma)$ is almost abelian of rank at most $\dim
\Sigma$; $\dim H^1(\Sigma;\Q) \leq \dim \Sigma$; and there is an
injection $H^*(\T^b;\Q) \hookrightarrow H^*(\Sigma;\Q)$ where $b=\dim
H^1(\Sigma;\Q)$. These constraints are ineffective for exotic tori.

In~\cite{RT}, Rudnev and Ten assume that a geodesic flow is completely
integrable with a non-degenerate first-integral map on an
$n$-dimensional compact manifold with first Betti number equal to
$n$. Non-degeneracy means, amongst other things, that the singular set
is stratified by the rank of the first integral map and each stratum
is a symplectic submanifold on which the system is completely
integrable. From these hypotheses, they deduce that there is a
lagrangian torus $F \subset T^* \Sigma$ such that the natural map
$\rho$ (figure \ref{fig:int1}) is a {\em homeomorphism}. Theorem 2 of
\cite{RT} states that $\rho$ is a diffeomorphism, but this is
mistaken. It is shown only that $\rho$ is a $1-1$ smooth map, hence by
invariance of domain, a homeomorphism. To prove that $\rho$ is a
diffeomorphism one must prove that the Maslov cocycle of $F$ vanishes,
or something equivalent. This is the first difficulty in proving
theorem \ref{thm:X}.

It should also be noted that either real-analyticity or non-degeneracy
is a very restrictive hypothesis on the first-integral map. In
\cite{Butler:2006a}, there is an example of a geometrically semisimple geodesic
flow on $T^*(\T^2 \times S^2)$ which is not completely integrable with
real-analytic (resp. non-degenerate) first integrals, nor is it
approximable by a real-analytically (resp. non-degenerately)
completely integrable system.

\subsubsection{Technical clarifications of theorem \ref{thm:X}}
One might also enquire if there is a convex hamiltonian $H$ which
enjoys an energy level $H^{-1}(c)$ which is geometrically
semisimple. If the sub-level $H^{-1}((-\infty,c])$ contains the zero
section of $T^* \Sigma$, then the answer is also {\em no}. Presumably,
the answer may change if the sub-level set does not contain the zero
section, but this is an open question (\cf \cite{RT}). The conditions
on the fibration $f: L \to B$ in the definition of geometric
semisimplicity may also be weakened somewhat and Theorem~\ref{thm:X}
continues to hold: one may require only that $\Gamma$ satisfy
condition (FI2) in Definition 9 of~\cite{Butler:2005} in place of
being a tamely-embedded polyhedron.

\subsubsection{ A sketch of the proofs }
If $g : X \times \R \to X$ is a flow, a point $x \in X$ is {\em
non-wandering} if, for any neighbourhood $U$ of $x$, $g_t(U) \cap U$
is non-empty for some $t > 1$. The set of non-wandering points for $g$
is denoted by $\nw{g}$ \cite{GPP}. It is proven that

\begin{thm}[\cf Theorem \ref{thm:van1}] \label{thm:int1}
Let $\Sigma$ be a smooth manifold and $H : T^*\Sigma \to \R$ a convex
hamiltonian with complete hamiltonian flow $\varphi$. If $F \subset
H^{-1}(c)$ is a lagrangian submanifold whose Maslov cocycle
vanishes and $\nw{\varphi|F}=F$, then $\rho$ (figure \ref{fig:int1})
is a smooth covering map. In particular, if $F$ is a torus, then
$\Sigma$ is finitely {\em smoothly} covered by a
torus.
\end{thm}

\begin{figure}[!h]
$
\xymatrix{
F\ \ar@{^{(}->}[rr]^{\iota_F={\rm incl.}} \ar@{->}[drr]_{\rho=\pi \cdot \iota_F} && \ T^*
\Sigma \ar@{->>}[d]^{\pi={\rm proj.}} \\ 
&& \Sigma. }
$
\caption{} \label{fig:int1} 
\end{figure}

This theorem is certainly known to experts, see~\cite{BP1992}
or \cite[section 2.5]{GPP} and references therein. What makes this
theorem crucial for the present note is that it provides a mechanism
whereby the smooth structure of $\Sigma$ enters: if, under the
hypotheses of Theorem~\ref{thm:X} one can prove that there must exist
a lagrangian standard $n$-torus $F \subset T^*\Sigma$ with vanishing
Maslov cocycle, and one can show that the degree of $\rho$ must be
$\pm 1$, then one has obtained a proof of the theorem. This is done in
sections
\ref{sec:nvI} and \ref{sec:X}. Section \ref{sec:mc} recalls the
definition and properties of the Maslov cocycle. Section \ref{sec:pe}
deals with projectively equivalent metrics on exotic tori and section
\ref{sec:Xhomo} proves Theorem~\ref{thm:KS} and related results.

\subsubsection*{Acknowledgements}
The author thanks V. Matveev, G. Paternain, A.  Ranicki, and
J. Robbins for their useful comments and discussion.

\section{The Maslov Cocycle} \label{sec:mc}
Let us recall the definition and construction of the Maslov cocycle as
it was introduced and developed in \cite{Ar1,Dui,Ar2}. An
interpretation of the Maslov cocycle as an obstruction class (a
primary difference) is also recalled.

\subsection{The Grassmannian of Lagrangian planes in $\C^n$}

Let $H$ be the standard hermitian inner product on $\C^n$: $H(z,w) =
\sum_{i=1}^n z_i \bar{w}_i$. This hermitian product is the sum of two real
quadratic forms on $\R^{2n} = \C^n$, the real symmetric part is a
euclidean inner product $g$ and the real skew-symmetric part is a
symplectic inner product $\omega$. The identity $H(z,w) = g(z,w) +
i\omega(z,w)$ shows that a real subspace $V \subset
\C^n$ is a real inner product space rel. $H$ iff $\omega|V = 0$. A
real inner-product subspace $V$ of $n$-dimensions is called a
lagrangian subspace; it is clear that $V$ is lagrangian iff there is a
basis $v_1,\ldots,v_n$ of $V$ such that $v = [ v_1 \cdots v_n ]$ is a
unitary matrix. This shows that the set of all lagrangian subspaces of
$\C^n$ is the homogeneous space $\Lambda_n = \U_n/\O_n$. It is a
standard exercise that every $2n$-dimensional symplectic vector space is
isomorphic to $(\R^{2n}, \omega)$.

The map $u \mapsto \det u^2=\exp(2\pi i \theta)$ induces a submersion
$\det^2 : \Lambda_n \to \U_1$. Let $\mu_o = \d \theta$, the standard
$\U_1$-invariant $1$-form on $\U_1$.
\begin{definition} \label{def:mc}
The {\em Maslov cocycle} $\mu$ is the pullback of
$\d \theta$ by $\det^2$. The Poincar\'e dual of $\mu$ is the {\em
Maslov cycle} $$\mc = \{ \lambda \in \Lambda_n\ :\ \lambda \cap i\R^n
\neq 0 \},$$ the set of lagrangian planes with a non-trivial
intersection with the plane $i\R^n$.
\end{definition}
The co-orientation of $\mc$ is defined by declaring that the closed curve
$c : [0,\pi] \to \Lambda_n$, $c : t \mapsto e^{it}\R \oplus i\R^{n-1}$
crosses $\mc$ positively at $t=\pi/2$. It is straightforward to see
that $\langle \mu, c \rangle = +1$, also.

\subsection{The bundle of lagrangian planes} \label{sec:lp}
If $(E,w) \to M$ is a symplectic vector bundle, then $E$ admits a
complex structure $J$ and a hermitian inner product $H$ such that
$(E_x,J_x,H_x,w_x)$ is isomorphic to $(\C^n,i,H,\omega)$ for all $x
\in M$. The associated bundle $\Lambda(E) \to M$ of lagrangian planes
is naturally defined.  Let $r,s : M \to \Lambda(E)$ be sections. The
primary obstruction to the existence of a homotopy between $r$ and $s$
is a cohomology class $d \in H^1(M;\{\pi_1(\Lambda_n)\})$, called the
primary difference. In general, the primary difference lies in a
cohomology group with twisted coefficients; because $B\U_n$ is simply
connected, the coefficients are untwisted in
$H^1(M;\{\pi_1(\Lambda_n)\})$. One may identify $\pi_1(\Lambda_n)$
with the integers by choosing the standard generator of
$\pi_1(\Lambda_n)$ to be the closed curve $c : [0,\pi] \to \Lambda_n$,
$c : t \mapsto e^{it}\R \oplus i\R^{n-1}$. With this convention, the
primary difference of $r$ and $s$ is a cohomology class $d\in
H^1(M;\Z)$.

When $M = \U_1$ and $E$ is a $\C^n$-vector bundle over $M$, then the
simple connectedness of $B\U_n$ implies that $E=M \times
\C^n$ and $\Lambda(E)=M\times \Lambda_n$. If $r,s : M \to \Lambda(E)$
are sections, then one has the cohomology classes $\mu_r=r^*(1 \times
\mu)$ and $\mu_s=s^*(1 \times \mu)$ and the primary difference equals
$$d = \mu_r - \mu_s.$$ Both $\mu_r$ and $\mu_s$ must be
multiples of the generator $\mu_o \in H^1(\U_1;\Z)$. It is clear that
this multiple is the degree of the composite maps 
$$\xymatrix{
U_1 \ar[r]^{r,s}\ar@/_3mm/[rr] &  \Lambda_n \ar@{->>}[r]^{\det^2} & U_1.
}
$$ 

By naturality, this characterizes the primary difference in all
cases: one can unnaturally trivialize the vector bundle $E \to M$ over
the $1$-skeleton of $M$ and apply the preceding to determine the
primary difference.

The Poincar\'e dual to the primary difference $d$ is a codim-$1$
cycle that is henceforth denoted by $\mc_\zeta$. By transversality,
one may assume that $r$ and $s$ are transversal sections $M \to
\Lambda(E)$; the set
$$\mc_\zeta = \left\{ m\in M\, :\, r(m) \cap s(m) \neq 0 \right\}$$ is
a smooth cycle that one may justifiably call a relative Maslov cycle.

The primary difference has a second, very important,
interpretation. Since $s$ and $J \circ s$ are homotopic, everywhere
transverse sections, $d$ is also the primary obstruction to $r$ being
homotopic to a section that is everywhere transverse to $s$.

\begin{remark} \label{rem:a}
1/ $(E,J,s)$ is a complex vector bundle with a real sub-bundle
$s$. Homotopy classes of sections of $\Lambda(E)$ therefore classify
the inequivalent real forms of $E$. These triples are classified, up
to isomorphism, by the set of homotopy classes of maps $[M, \U_n/{\rm
    O}_n]$ modulo the image of $[\Sigma^1 M,B\U_n]$ under a connecting
map~\cite[pp.s 97--99]{ML}. 2/ We have used the fact that a symplectic
vector bundle $E$ admits a non-natural complex structure $J$. The set
of such complex structures is contractible, so the non-naturality does
not affect homotopy invariants. 3/ In \cite{CGIP}, and in the preprint
of the present paper, one finds the statement that the existence of a
section of $\Lambda(E)$ implies its triviality. As pointed out by
Leonardo Macarini, this is incorrect. Indeed, here is a simple example
which shows that $\Lambda(E)$ may be non-trivial and have
sections. According to \cite[pp.s 97--99]{ML}, if $\Lambda(E)$ is
trivial, then $E$ is a trivial complex vector bundle. In our example,
$\Lambda(E)$ has a section but $E$ is non-trivial, so $\Lambda(E)$
must be non-trivial. Let $f : S^3 \to \SO_3$ be the canonical
projection map, let $j : \SO_3\to\U_3$ be the natural embedding. Let
$E\to S^4$ be the $\C^3$-vector bundle over $S^4$ whose clutching
function is $jf$. $E$ is a non-trivial $\C^3$ bundle since $jf
\in\pi_3(\U_3)$ is not trivial. Of course, $E$ admits a global
lagrangian sub-bundle, so $\Lambda(E)$ is non-trivial.

\end{remark}

\subsection{Cotangent bundles}

Let us specialize the constructions above. $E := T(T^*\Sigma)$ is a
symplectic vector bundle over $T^*\Sigma$ with the canonical
symplectic form $w = {\rm d}p \wedge {\rm d}q$. The footpoint
projection $\pi : T^*\Sigma \to \Sigma$ induces the {\em vertical
sub-bundle} $V = \ker({\rm d}\pi)$ of $E$. The fibres of $V$ are
lagrangian planes and the map $s(\theta)=V_{\theta}$ is a section of
the lagrangian grassmannian bundle $\Lambda(E)
\stackrel{\Pi}{\longrightarrow} T^*\Sigma$.

\begin{wrapfigure}[9]{r}[0pt]{0pt}
$
\xymatrix{
F\ \ar@{^{(}->}[r]^{\iota_F} & \ M\ \ar@{->>}[d]^f \ar@{^{(}->}[r]^{\iota_M} & T^* \Sigma\\
& B.
}
$
\caption{}\label{fig:4}
\end{wrapfigure} 
Let $M$ be a submanifold of $T^*\Sigma$ that is fibred by compact
lagrangian submanifolds, so that Figure \ref{fig:4} obtains where
$\iota_{\bullet}$ is an inclusion map, $f$ is a fibre-bundle map whose
fibres are lagrangian submanifolds and $F$ is a typical fibre.  This
hypothesis includes the possibility that $M$ itself is a lagrangian
submanifold. Since the fibres of $f$ are lagrangian submanifolds, for
each $\theta \in M$ the tangent space to the fibre of $f$ at $\theta$
is a lagrangian plane in $E_{\theta}$. Let $E_M = E|M$ and define a
section $r : M \to \Lambda(E_M)$ by
$$r(\theta) = \ker {\rm d}_{\theta}f, \qquad \forall\, \theta \in M.$$
The discussion of the relative Maslov cocycle and cycle from the
previous section applies to the present construction: one can compute
the primary difference $d$ between the section $r$ and the vertical
section $s|M$ of $\Lambda(E_M)$. As above, the relative Maslov cycle
$\mc_{\zeta}$ is the set of points where $\ker {\rm d}f$ and $\ker
{\rm d}\pi$ have a non-trivial intersection. If $F_{\theta}$ denotes
the fibre of $f$ through $\theta$, then
$$\mc_{\zeta} = \{ \theta \in M\ :\ {\rm rank\,} {\rm d}_{\theta}
\pi|F_{\theta} < \dim \Sigma
\}.$$
Since $F_{\theta}$ and $\Sigma$ are both $n$-dimensional manifolds,
$\mc_{\zeta}$ is the set of $\theta$ where $\pi|F_{\theta} :
F_{\theta} \to \Sigma$ fails to be a local diffeomorphism.

\section{A non-vanishing Maslov cocycle} \label{sec:nvI}

This section continues with the notation of the previous. Let us state
the main result of this section. Recall that $\nw{ \bullet }$ is the
non-wandering set of the flow $\bullet$.

\begin{thm} \label{thm:van1}
If $\mc_{\zeta}\, \cap\, \nw{\varphi|F} \neq \emptyset$, then $\mc_{\zeta}
\cap F$ is a non-torsion codimension-1 cycle on $F$.
\end{thm}

Note that this theorem requires only that $\mc_{\zeta}$ intersect the
chain recurrent set of $\varphi|F$, similar to \cite{BP1992}. However,
since $F$ will generally be a Liouville torus for this paper,
$\varphi|F$ will satisfy the somewhat stronger hypothesis here. The
proof has been included for the sake of completeness. The basic
underlying fact is that solution curves of convex hamiltonian systems
cross the Maslov cycle positively, as pointed out by
Duistermaat~\cite{Dui}.

\begin{proof}
Let $\theta \in \mc_{\zeta} \cap
  \nw{\varphi|F}$. By the convexity of $H$, there is an $s > 0$ such
  that
$$t \in [-s,s] \quad \textrm{and} \quad \varphi_t(\theta) \in \mc_{\zeta} \qquad
\implies \qquad t=0.$$
Since the non-wandering set $\nw{\varphi|F}$ is invariant, the points
$\theta^{\pm} = \varphi_{\pm s}(\theta)$ are non-wandering. Let
$U^{\pm}$ be neighbourhoods of $\theta^{\pm}$ that are disjoint from
$\mc_{\zeta}$. Since the points are non-wandering and on the same
orbit, there is a point $\theta' \in U^+$ and a $T>1$ such that
$\varphi_T(\theta') \in U^-$. 

Let $\gamma$ be the curve in $F$ obtained by concatenating the orbit
segment $\varphi_t(\theta) : t \in [-s,s]$, followed by an arc in
$U^+$ joining $\theta^+$ to $\theta'$, followed by the orbit segment
$\varphi_t(\theta') : t \in [0,T]$, followed by a segment joining
$\varphi_T(\theta')$ to $\theta^-$ in $U^-$. Since the segments of
$\gamma$ in $U^{\pm}$ are disjoint from $\mc_{\zeta}$, and the
remaining segments are $\varphi$-orbit segments, convexity implies
that
$$\#(\gamma, \mc_{\zeta}) \geq 0.$$
Since $\theta \in \gamma \cap \mc_{\zeta}$, the intersection number is
positive. 
\end{proof}

\begin{remark} \label{rem:b}
 If $F$ is a Liouville torus of a completely
integrable convex hamiltonian, then the Liouville-Arnold theorem
implies that $\nw{\varphi|F}=F$. Therefore, Theorem~\ref{thm:van1}
implies that the Maslov cocycle $\iota_F^*(d)$, if non-zero,
represents a non-torsion cohomology class in $H^1(F)$.
\end{remark}

\begin{corollary} \label{cor:van1}
Assume $\nw{\varphi|F}=F$. Then $\mc_{\zeta} \cap F$ is a trivial
cycle iff $\iota_F^*(d)$ is a trivial cocycle iff the map $\rho=\pi
\circ \iota_F$ in Figure \ref{fig:int1} is a local diffeomorphism.
\end{corollary}

\section{Exotic Tori and...} \label{sec:X}

\subsection{Geometric semisimplicity } \label{ssec:Xgs} 
A {\em topological} $n$-torus is a smooth manifold that is
homeomorphic to the {\em standard} $n$-torus $\T^n = \R^n/\Z^n$.  An
exotic $n$-torus is a topological $n$-torus that is not diffeomorphic
to $\T^n$.

\begin{prop} \label{prop:x1}
Let $\Sigma$ be a topological $n$-torus. If $H : T^*\Sigma \to \R$ is
a finitely geometrically semisimple convex hamiltonian, then there is
a lagrangian torus $F \subset T^*\Sigma$ such that the map $\rho$
$$
\xymatrix{
F\ \ar@{^{(}->}[rr]^{\iota_F={\rm incl.}} \ar@{->}[drr]_{\rho=\pi \cdot \iota_F} && \ T^*
\Sigma \ar@{->>}[d]^{\pi={\rm proj.}} \\ 
&& \Sigma. }
\eqno(\textrm{Figure \ref{fig:int1}})$$ has a non-zero degree.
\end{prop}

\begin{proof}
Let $f : L \to B$ be a lagrangian fibration, invariant under the
hamiltonian flow of $H$, such that $T^* \Sigma$ is the disjoint union
of $L$ and a closed, nowhere dense, tamely-embedded polyhedral
singular set $\Gamma$. Define the natural map $\xi$ by
\begin{equation} \label{eq:xi}
\xymatrix{L \ar@<-1mm>@/_3mm/[rrrr]_{\xi=\pi \cdot \iota_L}
\ar@{^{(}->}[rr]^{\iota_L=\textrm{incl.}} && T^* \Sigma \ar@{->>}[rr]^{\pi} &&
\Sigma}.
\end{equation}
By \cite[Lemma 15]{Butler:2005}, and the fact that both $F$ and
$\Sigma$ are topological $n$-tori, there is a component $L_i$ of $L$
such that $(\xi \cdot \iota_{L_i})_* : \pi_1(L_i) \to \pi_1(\Sigma)$
is almost surjective. Let us drop the subscript $i$ in the following
discussion; equivalently, let us assume that $L_i=L$.

By hypothesis, $\pi_1(B)$ is finite. The homotopy long exact sequence
for the fibration $(f,L,B)$ yields
$$
\xymatrix{
\cdots \ar[r] & \pi_2(B) \ar[r]^{\partial_*} & \pi_1(F) \ar[r]^{\iota_{F,L,*}} & \pi_1(L) \ar[r]
& \pi_1(B) \ar[r]^{f_*} & 1,
}
$$ so $\pi_1(L)$ contains the finite-index subgroup
$\pi_1(F)/\partial_* \pi_2(B)$. One concludes that, since
$\rho=\xi|F$, the image of $\pi_1(F)$ under the map $\rho_* : \pi_1(F)
\to \pi_1(\Sigma)$ is a finite index subgroup. Since both $F$ and
$\Sigma$ are topological $n$-tori, $\rho$ has a non-zero degree.
\end{proof}

We continue with the hypotheses of Proposition~\ref{prop:x1}.

\begin{prop} \label{prop:x2}
If $\deg \rho \neq 0$ (Figure \ref{fig:int1}), then $\rho$ is a local
diffeomorphism.
\end{prop}

\noindent
(It is known that there is a unique PL structure on the topological
$n$-torus, a fact that is used without further reference.)

\begin{proof}
(\cf \cite{Viterbo}) It is known that the smooth structure of the
  topological $n$-torus $\Sigma$ is determined by a unique cohomology
  class $\sigma$ contained in the cohomology group $\bigoplus_{i \leq
    n} H^i(\Sigma; \Gamma_i)$, where $\Gamma_i = \pi_{i}(PL/O)$ is the
  group of homotopy classes of maps from $S^{i}$ into the classifying
  space of stable $PL$-structures modulo smooth structures;
  equivalently, $\Gamma_i$ is the group of smooth structures on the
  topological $i$-sphere for $i\geq 7$ and $0$ for $i<7$
  \cite[p. 236]{Wall}. This correspondence is natural with respect to
  local diffeomorphisms, so if $p : \Sigma' \to \Sigma$ is a local
  diffeomorphism and $\Sigma$ (resp. $\Sigma'$) is a topological
  $n$-torus whose smooth structure is determined by the cohomology
  class $\sigma$ (resp. $\sigma'$), then $p^*\sigma = \sigma'$.

Since the cohomology class $\sigma$ lies in a finite group, it has
finite order. Therefore, if $p : \Sigma' \to \Sigma$ is a finite
covering whose degree divides the order of $\sigma$, then $\sigma' =
p^* \sigma$ must vanish. Thus $\Sigma'$ is diffeomorphic to the
standard $n$-torus. It is clear that such coverings $p$ exist, so let
us choose one such covering.

\begin{center}
\begin{figure}[htb]
$
\xymatrix{
F'\ 
\ar @{->} @/^5mm/[rrrrrr]^{\phi}
\ar@{^{(}->}[rr]_{\iota_{F'}={\rm incl.}} \ar@{->}[drr]_{\rho'=\pi' \cdot \iota_{F'}} && \ T^*
\Sigma' \ar@{->>}[d]^{\pi'={\rm proj.}} \ar@{->>}[rr]^{P={\rm proj.}} && T^*
\Sigma \ar@{<-^{)}}[rr]_{\iota_{F}={\rm incl.}} \ar@{->>}[d]_{\pi={\rm proj.}} && F \ar@{->}[dll]^{\rho=\pi \cdot \iota_F} \\ 
&& \Sigma' \ar@{->>}[rr]^{p} && \Sigma. }
$
\caption{The pullback diagram of Figure \ref{fig:int1}.} \label{fig:int1+}
\end{figure}
\end{center}

From Figure \ref{fig:int1}, one gets the pullback diagram where $F'$
is a connected component of $P^{-1}(F)$ and $\phi$ is the covering map
induced by $P$, see Figure \ref{fig:int1+}. Because $F$ is
diffeomorphic to the standard $n$-torus and $F'$ is a finite covering
of $F$, $F'$ is also diffeomorphic to the standard $n$-torus--this
follows from the above-mentioned classification of smooth structures
and their naturality under coverings.

Since $\deg \rho$ is non-zero by hypothesis and $\deg \rho' \cdot \deg
p = \deg \phi \cdot \deg \rho$, the degree of $\rho'$ is non-zero. The
map $\iota_{F'}$ is an embedding by naturality. Therefore,
$\iota_{F'}$ is a lagrangian embedding of the standard $n$-torus $F'$
into the cotangent bundle of the standard $n$-torus $T^* \Sigma'$ such
that the induced map $\rho'$ has a non-zero
degree. Viterbo~\cite[Corollary 3]{Viterbo} proved that in this case the Maslov
cocycle $\iota_{F'}^*(d')$ is cohomologically trivial. By the remark
following Theorem~\ref{thm:van1}, the Maslov cycle $\mc_{\zeta'} \cap
F'$ is therefore empty. Thus $\rho'$ is a local diffeomorphism. The
commutativity of Figure \ref{fig:int1+} shows that this forces $\rho$
to be a local diffeomorphism.
\end{proof}

\begin{prop} \label{prop:x3}
The degree of $\rho$ is $\pm 1$. Hence $\rho$ is a diffeomorphism of
the standard $n$-torus $F$ with $\Sigma$.
\end{prop}

\begin{proof}
Let $\lambda$ (``$={\rm p}\cdot\d {\rm q}$'') be the Liouville
$1$-form of $T^*\Sigma$ and define
\begin{equation} \label{eq:a}
\alpha = \iota_F^*(\lambda).
\end{equation}
Since $F$ is a lagrangian manifold, $\alpha$ is a closed $1$-form on
$F$. 

Let $\Deck$ be the deck transformation group of the covering map
$\rho$. Let $\forms{k}{F}$ be the vector space of smooth $k$-forms on
$F$. $\Deck$ acts linearly on $\forms{k}{F}$ via pullback; let
$\forms{k}{F}^\Deck$ denote the fixed-point set of $\Deck$'s action on
$\forms{k}{F}$. It is a well-known fact that
$$\frac{\ker \d | \forms{k}{F}^\Deck}{ \im \d|\forms{k+1}{F}^\Deck } \stackrel{\rho^*}{\cong}
H^k_{\textrm{de Rham}}(F/\Deck) = H^k_{\textrm{de Rham}}(\Sigma).$$

Since $F$ and $\Sigma$ are topological $n$-tori, their cohomology
groups are isomorphic. From these facts, there is a decomposition
\begin{equation} \label{eq:a1}
\alpha = \alpha_0 + \alpha_1
\end{equation}
where $\alpha_0 \in \forms{1}{F}^\Deck$ is cohomologous to $\alpha$
and therefore $\alpha_1 = \d h$ is exact.

Observe that for all $x \in F$ and $\gamma\in \Deck$
$$\alpha_x - \gamma^*\alpha_{\gamma(x)} = \d h^\gamma_x$$ where
$h^\gamma(x) = h(x)-h(\gamma(x))$ is a smooth function. Since $F$ is
compact, there is a critical point $x=x_{\gamma}$ of $h^\gamma$, so
that
\begin{equation} \label{eq:aav}
\alpha_x - \gamma^*\alpha_{\gamma(x)} = 0.
\end{equation}

To complete the proof, it is claimed that equation (\ref{eq:aav})
implies that $\gamma = 1$. This proves that the deck transformation
group $\Deck$ is trivial, whence $\rho$ is a diffeomorphism.

To prove the claim, recall that the Liouville $1$-form $\lambda$ at
$\theta \in T^*\Sigma$ is equal to the $1$-form $\left(\d_{\theta}\pi\right)^*
\theta \in T^*_{\theta}(T^*\Sigma)$. Therefore, for each $x\in F$,
\begin{equation} \label{eq:ax}
\alpha_x = \left( \d_x \iota_F \right)^* \cdot \left( \d_{\iota_F(x)}
\pi \right)^* (\iota_F(x)) = \left( \d_x \rho \right)^*(x),
\end{equation}
where in the second step the identity $\rho = \pi \circ \iota_F$ has
been used and the innocuous inclusion map dropped. Equation
(\ref{eq:ax}) implies that for all $\gamma \in \Deck$
\begin{equation} \label{eq:acx}
\gamma^* \alpha_{\gamma(x)} = \left( \d_x \gamma \right)^* \cdot
\left(  \d_{\gamma(x)} \rho  \right)^* (\gamma(x)) =
\left( \d_x \rho\right)^*( \gamma(x) ),
\end{equation}
where $\rho \circ \gamma = \rho$ and the fact that $x,\gamma(x) \in
T^*_{\rho(x)}\Sigma$ has been used. Therefore
\begin{equation} \label{eq:af}
\alpha_x - \gamma^*\alpha_{\gamma(x)} = \left( \d_x \rho  \right)^*
\left( x - \gamma(x)  \right).
\end{equation}
Since $\rho$ is a local diffeomorphism, equation (\ref{eq:af}) shows
that $\alpha_x-\gamma^*\alpha_{\gamma(x)} = 0$ iff
$x=\gamma(x)$. Since $\Deck$ acts freely on $F$, equation
(\ref{eq:aav}) therefore implies that all elements of $\Deck$ are $1$.
\end{proof}

\begin{proof}[Theorem~\ref{thm:X}]
Proposition~\ref{prop:x3} proves Theorem~\ref{thm:X}.
\end{proof}

\begin{remark} \label{rem:c}
Lalonde and Sikorav~\cite[p. 19]{LS} ask if the map $\rho$ in
figure~\ref{fig:int1} has $\deg \rho = \pm 1$ or possibly just $\neq
0$, when $F$ is an exact lagrangian submanifold (\ie when the $1$-form
$\alpha$ in equation (\ref{eq:ax}) is exact). They prove that $\deg
\rho=\pm 1$ and the Maslov class $\iota_F^*(d)$ vanishes when
$F=\Sigma=\T^n$. In a similar vein, Bialy and Polterovich~\cite{BP1}
prove that an invariant lagrangian $2$-torus $F$ contained in the unit
co-sphere bundle of the $2$-torus has $\deg \rho=\pm 1$ iff $F$ is the
disjoint union of lifts of globally minimizing unit-speed geodesics on
the $2$-torus. These results are sharpened in \cite{BP2}. In
\cite{P1991}, Polterovich proved that if $F$ is an exact lagrangian
torus in $T^*\T^n$, then $F$ is a graph of a closed $1$-form; see also
\cite[Corollary II]{BP1992} and \cite[Theorem
  1.1]{BP1992.1}. Propositions~\ref{prop:x2}--\ref{prop:x3} may be
viewed as an extension of each of these results.
\end{remark}

\section{Projectively Equivalent Metrics} \label{sec:pe} 

\subsubsection{Preamble } \label{sssec:intro} 
Let $\Sigma$ be a smooth $n$-dimensional manifold and let $\g, \bg$ be
smooth riemannian metrics on $\Sigma$. These metrics are said to be
projectively equivalent if their geodesics coincide as unparameterized
curves. Projective equivalence is related to complete integrability in
the following manner.

Define a $\bg$-self-adjoint $(1,1)$ tensor field $\L$ by
\begin{equation} \label{eq:L}
\L = \left( \frac{\det(\bg)}{\det(\g)}  \right)^{\frac{1}{n+1}} \times
\bg^{-1} \cdot \g, 
\end{equation}
where one views the metrics as self-adjoint bundle isomorphisms
$T\Sigma \to T^*\Sigma$. At each point $x \in \Sigma$, $\L$ has $n$ real
eigenvalues, and one can define continuous functions $\lambda_i$ by
declaring $\lambda_i(x)$ to be the $i$-th largest eigenvalue of $\L$
at $x$. The metrics are said to be strictly non-proportional at $x$ if
$\L$ has $n$ distinct eigenvalues there. The functions $\lambda_i$ are
smooth in a neighbourhood of such an $x$.

Define a polynomial family of $(1,1)$ tensor fields by
\begin{equation} \label{eq:St}
\S_\tau = \adj{\L-\tau},
\end{equation}
where $\adj{\bullet}$ is the classical adjoint matrix and $\tau$ is a
real number. From these tensor fields, one obtains functions
\begin{equation} \label{eq:It}
I_\tau(x,v) = \left\langle \g \cdot \S_\tau \cdot v , v \right\rangle,
\hspace{10mm} \forall (x,v)\in T\Sigma.
\end{equation}
Let $J_\tau = I_\tau \cdot \g^{-1}$ be the pullback of these functions
to $T^*\Sigma$. Note that the lagrangian of the riemannian metric of $\bg$
is $I_0$, while that of $\g$ equals $\lim_{\tau\to\infty} \tau^{-n+1}
I_\tau$.

\begin{qthm}[Topalov-Matveev 1998] \label{thm:MT}
The family $\left\{ J_\tau   \right\}_{\tau\in\R}$ is a Poisson
commuting family. If there exists a point $x\in M$ where $\L$ has
$n$ distinct eigenvalues, then the geodesic flow of $\g$ is completely
integrable.
\end{qthm}

This theorem, along with the theorem of Levi-Civita which establishes
a normal form for the metrics in the neighbourhood of a regular point,
suffice to prove the following:

\begin{thm} \label{thm:nope}
If $\Sigma$ is a topological $n$-torus and $\g,\bg$ are projectively
equivalent metrics that are strictly non-proportional at a single
point, then $\Sigma$ is diffeomorphic to the standard $n$-torus.
\end{thm}

\begin{proof}
From the discussion in Proposition \ref{prop:x2}, there is a finite
covering $p : \Sigma' \to \Sigma$ where $\Sigma'$ is diffeomorphic to
the standard $n$-torus. The metrics $p^*\g, p^* \bg$ are also
projectively equivalent and strictly non-proportional at some
point. 

Say that $\g_m,\bg_m$ are `model' metrics on the standard torus $\T^n
= \R/\Z \times \cdots \times\R/\Z$ if 
\begin{equation} \label{eq:model}
\begin{split}
\g_m &= \sum_{i=1}^n \Pi_i\, \d x_i^2,\\
\bg_m &= \sum_{i=1}^n \rho_i \Pi_i\, \d x_i^2, \hspace{10mm}
\textrm{where}\\
\Pi_i &= (-1)^{n-i-1}\, \prod_{j\neq i} (\lambda_i-\lambda_j),\hspace{5mm}
\rho_i^{-1} = \lambda_i \cdot \lambda_1 \cdots \lambda_n,
\end{split}
\end{equation}
and $\lambda_i=\lambda_i(x_i)$ is a function of the $i$-th coordinate
alone and
\begin{equation} \label{eq:ineq}
i<j \hspace{5mm} \implies \hspace{5mm} \lambda_i(x) < \lambda_j(y)
\hspace{5mm}\forall x,y.
\end{equation}
By \cite[Theorem 7]{Matveev}, there is a diffeomorphism $h : \Sigma'
\to \T^n$ which is an isometry of $p^*\g,p^*\bg$ with model metrics
$\g_m,\bg_m$. Henceforth, it is assumed without loss of generality
that $h$ is the identity, $\Sigma'=\T^n$ and $p^*\g=\g_m,
p^*\bg=\bg_m$.

In the coordinate system on $\Sigma'$, one computes that 
\begin{equation} \label{eq:S}
\begin{split}
\S_\tau &= \sum_{i=1}^n \mu_i(\tau)\, \frac{\partial \ }{\partial x_i}
\otimes \d x_i, \hspace{10mm} \mu_i(\tau) = \mu_i(\tau;x) = \prod_{j\neq i}
(\lambda_j-\tau),\\
J_\tau &= \sum_{i=1}^n \mu_i(\tau)\, \Pi_i^{-1}\, y_i^2,
\end{split}
\end{equation}
where $(x_i,y_i)$ are canonical coordinates on $T^*\Sigma'$. The
hamiltonian vector field $X_{J_\tau}$ equals
\begin{equation} \label{eq:xj}
X_{J_\tau} = \left\{ 
\begin{array}{lcl}
\dot{x}_i &=& \frac{2\mu_i y_i}{\Pi_i},\\ \dot{y}_i &=& \lambda_i' \
\sum_{j\neq i} \left(
\frac{\lambda_j-\tau}{(\lambda_i-\tau)(\lambda_j-\lambda_i)} \right)\cdot
\frac{\mu_jy_j^2}{\Pi_j}.
\end{array}
\right.
\end{equation}
Define $\H_{(x,y)} \subset T_x\Sigma'$ to be the subspace spanned by the
projection of the tangent vectors $X_{J_\tau}(x,y)$ to the base. If one
chooses real numbers $t_1<\cdots<t_n$, then it is clear that
\begin{equation} \label{eq:h}
0 \neq \det \left[ 
\begin{array}{ccc}
\mu_1(t_1) & \cdots & \mu_n(t_1)\\
\vdots & \ddots & \vdots\\
\mu_1(t_n) & \cdots & \mu_n(t_n)
\end{array}
\right] 
\times \prod_{i=1}^n
\frac{y_i}{\Pi_i} \hspace{5mm} \implies \hspace{5mm} \dim \H_{(x,y)}=n.
\end{equation}
Equation (\ref{eq:ineq}) implies that for all $x\in\Sigma'$, the
eigenvalues $\lambda_i$ are pairwise distinct, so the theory of
Lagrange interpolation shows that the polynomials $\mu_i(\tau)$ are
linearly independent, whence the determinant on the left is nowhere
zero. To complete the proof of the theorem, it therefore suffices to
show that there is a lagrangian torus $F' \subset T^* \Sigma'$ such that
$y_1 \cdots y_n$ does not vanish on $F'$.

\medskip
\noindent{\em The image of the first-integral map}. Let $\V$ be the vector
space of polynomials in $\tau$ of degree at most $n-1$. The map $(x,y)
\mapsto J_\tau(x,y)$ defines a smooth map $J : T^*\Sigma' \to \V$. The
following claim is essential to describe the image of $J$.

\begin{claim} \label{claim:a}
Let $x \in \Sigma'$ and $t_i=\lambda_i(x)$. If
$\tau_i \in [t_i,t_{i+1}]$ for all $i$, then there exists $a_i \geq
0$ with $\sum_{i=1}^n a_i = 1$, such that
\begin{equation} \label{eq:p}
p(\tau) = \sum_{i=1}^n a_i \mu_i(\tau;x)
\end{equation}
vanishes at $\tau_i$ for all $i$. If $\tau_i \in (t_i,t_{i+1})$ for
 all $i$, then $a_i>0$ for all $i$.
\end{claim}

\begin{chec}
Let $\tau_1 \leq \tau_2 \leq \cdots \leq
\tau_{n-1}$ be real numbers and define
\begin{equation} \label{eq:sigma}
\sigma(\tau) := \prod_{i=1}^{n-1} (\tau_i-\tau).
\end{equation}
A computation shows that if $p(\tau)$ in equation (\ref{eq:p})
vanishes at $\tau_i$ for all $i$, then the coefficients $a_i$ may be
written as
\begin{equation} \label{eq:ai}
a_i = \frac{\sigma(t_i)}{\mu_i(t_i)} \hspace{5mm}\forall
i=1,\ldots,n,
\end{equation}
where $\sum a_i=1$. Comparison of equations
(\ref{eq:sigma},\ref{eq:S}) shows that $a_i \geq 0$ (resp. $a_i > 0$)
for all $i$ iff $\tau_i \in [t_i,t_{i+1}]$ (resp. $\tau_i
\in (t_i,t_{i+1})$) for all $i$.
\end{chec}

\medskip
To continue the mainline of the proof: Let $\ul_i$ (resp. $\bl_i$) be
the maximum (resp. minimum) value attained by $\lambda_i$. Since the
eigenvalues $\lambda_i$ are everywhere distinct by (\ref{eq:ineq}),
there exists $\tau_i \in (\ul_i,\bl_{i+1})$ for all $i$. Since
$\lambda_i(x) \leq \ul_i < \tau_i < \bl_{i+1} \leq \lambda_{i+1}(x)$
for all $x\in\Sigma'$ and $i$, the claim establishes that for all
$x\in\Sigma'$, there is a $y\in T^*_x\Sigma'$ such that $J_{\tau}(x,y)
= p(\tau)$ where $p$ has roots $\tau_1, \ldots, \tau_{n-1}$. Moreover,
the coefficients $a_i = y_i^2/\Pi_i$ are everywhere non-zero by the
same claim. From (\ref{eq:h}), one sees that the canonical projection
map $\rho' : F' \to \Sigma'$ is a local diffeomorphism, where $F' =
J^{-1}(p)$. Proposition \ref{prop:x3} implies that $\rho'$ is a
diffeomorphism. Therefore, the map $\rho$ in figure \ref{fig:int1+} is
a diffeomorphism, and $\Sigma$ is therefore diffeomorphic to the
standard $n$-torus.
\end{proof}

\begin{remark} \label{remark:ge}
Claim \ref{claim:a} implies that the image of the first-integral map
$J$ contains the polynomials of the form
\begin{equation} \label{equation:p1}
q(\tau) = a \times \prod_{i=1}^{n-1} (\tau_i-\tau), \hspace{5mm} a>0,\
\tau_i\in (\bl_i,\ul_{i+1}).
\end{equation}
A connected component of the pre-image $J^{-1}(q)$ of such a
polynomial is a regular lagrangian torus whose projection to the base
$\Sigma'$ is a diffeomorphism. If one of the roots $\tau_i$ of $q$ lie
in $[\bl_i,\ul_i]$, and $q$ is a regular value of $J$, then the Maslov
cycle of each component of $J^{-1}(q)$ is non-trivial and one can see
from the claim that the projection of $J^{-1}(q)$ does not cover the
base.

Additionally, the claim shows that if $\tau_i \in [\ul_i,\bl_{i+1}]$
for all $i=1,\ldots,n-1$, then there exists some $(x,y)\in T^*\Sigma'$
such that $J_\tau(x,y) = q(\tau)$. This shows that the image of $J$ is
the set of all polynomials
\begin{equation} \label{equation:p2}
\begin{split}
q(\tau) &= a \times \prod_{i=1}^{n-1} (\tau_i-\tau),
\hspace{10mm} {\rm such\ that\ }a\geq 0,  \hspace{2mm}
\tau_i \in [\ul_i,\bl_{i+1}]\\
&{\rm and\ } \tau_1 \leq \tau_2 \leq
\cdots \leq \tau_{n-1}.
\end{split}
\end{equation}

\medskip
In \cite{Matveev}, Matveev proves Theorem~\ref{thm:nope} in a
different manner. The model metric's Levi-Civita coordinate system on
$\Sigma'$ (see \eqref{eq:model}) induces an integral affine structure
on $\Sigma$ whose holonomy group lies in the orthogonal
group. Therefore $\Sigma$ admits a flat riemannian metric, so the
second Bieberbach theorem implies that $\Sigma$ is diffeomorphic to
the standard $n$-torus. Finally, \cite[Section 9]{MT} contains the
proof of a theorem similar to Theorem~\ref{thm:nope}, also.
\end{remark}

\section{Exotic smooth $7$-manifolds} \label{sec:Xhomo}
\def\g{{\mathfrak g}}

\subsection{Witten-Kreck-Stolz manifolds } \label{ssec:KS} 
Let $k,l$ be coprime integers that are both non-zero. The action of
$\U_1$ on $S^5 \times S^3$ by
\begin{equation} \label{eq:tor}
\forall z\in S^1, x\in S^5, y\in S^3: \hspace{5mm} z\cdot(x,y) =
(z^k\cdot x, z^l \cdot y)
\end{equation} 
is free. Let $M_{k,l}$ be the orbit space
of this action; it is a compact simply connected
$7$-manifold. Equivalently, let $G=\U_3 \times \U_2$ and let $U\subset
G$ be the subgroup isomorphic to $\U_2 \times \U_1^2$ defined by
\begin{equation} \label{eq:U}
U = \left\{ z^k \cdot \left[ \begin{array}{cc}
1 & 0\\
0 & a
\end{array}
\right] \oplus z^l \cdot \left[ 
\begin{array}{cc}
1 & 0\\
0 & w
\end{array}
\right]\ :\   a\in\U_2,\ z,w\in\U_1  \right\}.
\end{equation}
The manifold $M_{k,l}$ is $G$-equivariantly diffeomorphic to the
homogeneous space $G/U$. 

These manifolds have been studied by Witten in the context of
Kaluza-Klein theory and by Wang and Ziller, who constructed Einstein
metrics on each manifold with positive scalar curvature
\cite{Witten,WZ}. In~\cite{KS2}, Kreck and Stolz classify the
manifolds up to homeomorphism and diffeomorphism. As a consequence of
the triviality of $H^3(M_{k,l};\Z_2)$, it is known from smoothing
theory that if $M_{k',l'}$ is homeomorphic to $M_{k,l}$, then the
former is diffeomorphic to the latter connect-summed with an exotic
$7$-sphere. There are, therefore, at most $28$ oriented diffeomorphism
classes within any homeomorphism class. Combining Corollary D and the
Remark preceding it in \cite{KS2}, we have

\begin{qthm}[Kreck-Stolz 1988] \label{thm:KS-2}
If $l=0\bmod 4$, $l=0,3,4 \bmod 7$, $l\neq 0$ and $(k,l)=1$, then
the homeomorphism class of $M_{k,l}$ has $28$ diffeomorphism
classes. Each diffeomorphism class is represented by an $M_{k',l'}$ for
suitable $k',l'$.
\end{qthm}

The simplest example satisfying the hypotheses of the theorem is the
manifold $M_{1,4}$. By theorem B of \cite{KS2}, $M_{k',l'}$ is
homeomorphic to $M_{1,4}$ iff $l'=\pm 4$ and $k'=1 \bmod 32$; on the
other hand, $M_{k',l'}$ is diffeomorphic to $M_{1,4}$ iff $l'=\pm 4$
and $k'=1 \bmod 28 \times 32$. Thus $M_{32t+1,4}$ enumerates all
diffeomorphism classes of $M_{1,4}$ for $t=0,\ldots,27$.

\medskip

On the other hand, Mykytyuk and Panasyuk have studied the
integrability of the canonical quadratic hamiltonian on homogeneous
spaces. Theorem 3.10 of \cite{MP} implies

\begin{qthm}[Mykytyuk-Panasyuk 2004] \label{thm:MP}
Let $G$ be a compact reductive Lie group. Let $K\subset G$ be the
stabilizer of some element $a\in \g^*$. If $U \subset K$ contains the
identity component of $[K,K]$, then the quadratic hamiltonian on
$T^*(G/U)$ induced by a bi-invariant metric on $G$ is completely
integrable with real-analytic integrals.
\end{qthm}

To apply this theorem to the Witten-Kreck-Stolz manifolds, let
$$K = \left\{ u \cdot \left[ 
\begin{array}{cc}
1 & 0\\
0 & a
\end{array}
\right] \oplus v \cdot \left[ 
\begin{array}{cc}
1 & 0\\
0 & w
\end{array}
\right]\ :\ u,v,w \in \U_1,\ \ a\in \SU_2  \right\},$$
which is easily seen to be the stabilizer subgroup under the coadjoint
action of the element
$$
\left[
\begin{array}{cc}
i & 0\\
0 & 0
\end{array}
\right] \oplus 
\left[ 
\begin{array}{cc}
i & 0\\
0 & -i
\end{array}
\right] \in \u_3^* \oplus \u_2^*.$$ Since the subgroup $U$, defined in
equation (\ref{eq:U}) above, contains $[K,K]$ and $U \subset K$, the
Mykytyuk-Panasyuk theorem is applicable to each Witten-Kreck-Stolz
manifold. This proves

\begin{thm} \label{thm:Ex}
Each homogeneous space $M_{k,l}$ has a real-analytically completely
integrable convex hamiltonian on its cotangent bundle. In particular,
if $k$ and $l$ satisfy the conditions of Theorem~\ref{thm:KS-2}, then
each diffeomorphism class of manifolds homeomorphic to $M_{k,l}$ has
such a real-analytically integrable convex hamiltonian.
\end{thm}

\begin{remark} \label{rem:d}
There is a more pedestrian approach to proving
Theorem~\ref{thm:Ex} which uses the above description of $M_{k,l}$ as
the quotient of $S^5 \times S^3$ by a subgroup of the torus $V \cong \U_1
\subset \U_1 \times \U_1$ (equation (\ref{eq:tor})). Introduce the notation
\begin{align} \label{eq:ts}
& \begin{array}{lcl}
T^*S^5 &=& \left\{ (\x,\y) \in \C^3 \times \C^3\ :\ \x^* \x = 1,
\re{\y^*\x} = 0  \right\},\\
T^*S^3 &=& \left\{ (\w,\z) \in \C^2 \times \C^2\ :\ \w^* \w = 1,
\re{\z^*\w} = 0  \right\},
\end{array}
\end{align}
where $\x^*$ denotes the hermitian transpose of the vector $\x$ (this
use of ${}^*$ is confined to the present remark). The momentum map of
the $G=\U_3 \times \U_2$ action on $T^*(S^5 \times S^3)$ equals, in
these coordinates
\begin{equation} \label{eq:mmU3U2}
\mm_G = \frac{1}{2} \times \left( \x\y^* - \y\x^* \right) \oplus
\frac{1}{2} \times \left( \w\z^* - \z\w^* \right), \hspace{5mm} \mm_G :
T^*(S^5 \times S^3) \to \u_3 \oplus \u_2,
\end{equation}
while the momentum map of the subgroup $V$ equals
\begin{equation} \label{eq:mmV}
\mm_V = k \times \y^*\x + l \times \z^*\w,
\end{equation}
where we identify the Lie algebra of $V$ with $i\R$. Define $\pi_{a,b}
: \u_a \to \u_b$ to be the orthogonal projection onto the subalgebra
$\u_b \subset \u_a$, which we embed in the standard way in the lower
right-hand corner. Define the following functions for $\xi\oplus\eta
\in \u_3\oplus \u_2$ by
\begin{align} \label{al:f}
& \left\{ 
\begin{array}{lclclcl}
f_1 &=& -i \cdot \pi_{3,1}(\xi), & \hspace{5mm} & f_2 &=& -i
\cdot \Trace{\pi_{3,2}(\xi)},\\
f_3 &=& -i \Trace{\xi}, & \hspace{5mm} & f_4 &=& \frac{1}{2}
\cdot \Trace{ \pi_{3,2}(\xi)^2 },\\
f_5 &=& \frac{1}{2}
\cdot \Trace{ \xi^2 }, & & f_6 &=& -i \cdot \pi_{2,1}(\eta),\\
f_7 &=& -i \Trace{\eta}, & &  f_8 &=& \frac{1}{2} \cdot \Trace{ \eta^2 },
\end{array}
\right.
\end{align}
and define to be $H_a = f_a \circ \mm_G$. The functions $H_a$ are all
in involution, a fact that is easily confirmed by the observation that
$\pi_{a,b}$ is a Poisson map (it is the transpose of the inclusion
$\iota_{a,b} : \u_b \to \u_a$), and the functions being pulled-back
are Casimirs. A simple computation shows that the eight hamiltonians
are functionally independent at $(\x_o,\y_o)=(2/9,1/9,2/9,i,-4i,i)$
and $(\w_o,\z_o)=(3/5,4/5,4i,-3i)$. One also notes that $H=H_5+H_8$ is
a convex hamiltonian that equals $-\frac{1}{4}\times \left( |\y|^2 +
|\z|^2 - (\y^*\x)^2 - (\z^*\w)^2 \right)$ and that
$H'=H_5+H_8+\frac{1}{4}\left( H_3^2 + H_7^2 \right) = -\frac{1}{4}
\times \left( |\y|^2 + |\z|^2 \right)$.

One observes that 
\begin{equation} \label{eq:mmV2}
\mm_V = i k \times H_3 + i l \times H_7,
\end{equation}
and that $T^*M_{k,l}$ is canonically symplectomorphic to
$\mm_V^{-1}(0)/V$. It is observed that the point
$(\x_o,\y_o,\w_o,\z_o)$ at which the eight hamiltonians are
functionally independent lies in $\mm_V^{-1}(0)$. Therefore, since
$kl\neq 0$, the seven hamiltonians $H_a$, $a\neq 3$, descend to
$T^*M_{k,l}$ to give seven functionally independent real-analytic
hamiltonians that are in involution. The convex hamiltonian $H=H_5 +
H_8$ is the desired real-analytically integrable convex hamiltonian on
$T^*M_{k,l}$.  This provides an alternative, more computational, proof
of Theorem~\ref{thm:Ex}.
\end{remark}

\subsection{A Lemma on integrability} \label{ssec:int} 
A second class of well-known $7$-manifolds are the Eschenburg spaces,
which include the well-known Aloff-Wallach spaces~\cite{Esch,AW}. In
the next subsection, these spaces are described more completely. This
subsection proves a Lemma about the existence of real-analytically
completely integrable convex hamiltonians on spaces constructed like
the Eschenburg spaces. This Lemma is new and possibly of independent
interest. 

Let $G$ be a compact Lie group and let $G \times G$ act on $G$ by the
convention $$\forall h=(h_1,h_2)\in G\times G, g \in G: \qquad h\cdot
g = h_1 g h_2^{-1}.$$ This can be understood as a left action on $G$
by the group $H=G_+ \times G_-$, where $G_{+/-}$ is $G$ equipped with
left/right multiplication respectively.

Let $U \subset G_+ \times G_-$ be a closed subgroup that acts freely
on $G$. There are naturally induced maps
$$
\begin{array}{lr}
\xymatrix{
&&  \ \ \ G_+=G_+ \times 1 \ar@<-1ex>@{_{(}->}[d]_{{\rm incl.}=\jj_+}\\
U\ \ \ar@{^{(}->}[rr]^{\iota={\rm incl.}} \ar@<1ex>[rru]^{\iota_+}
\ar@<-1ex>[rrd]^{\iota_-} && G_+ \times G_- = H\ar@{->>}[u]_{\pi_{+}={\rm proj.}} \ar@{->>}[d]^{\pi_-={\rm proj.}}\\
&& \ \ \ G_-=1 \times G_- \ar@<1ex>@{^{(}->}[u]^{{\rm incl.}=\jj_-}
}
&
\xymatrix{
&&  \ \ \ \g_+^*=\g_+^* \times 0 \ar@<-1ex>[lld]_{\iota_+^*} \ar@<-1ex>@{<<-}[d]_{\jj_+^*}\\
\u^*\ \ \ar@{<<-}[rr]^{\iota^*={\rm incl.}^*}  && \g_+^* \times \g_-^*
=\h^*\ar@{<-_{)}}[u]_{\pi_{+}^*} \ar@{<-^{)}}[d]^{\pi_-^*}\\
&& \ \ \ \g_-^*=0 \times \g_-^*, \ar@<1ex>[llu]^{\iota_-^*} \ar@<1ex>@{<<-}[u]^{\jj_-^*}
}
\end{array}
$$ where notation is abused and the induced map on Lie algebras is
denoted by the same symbol. 

For a subgroup $S$ of $H$, let $\mm_S : T^* G \to {\mathfrak s}^*$ be
the momentum map of $S$'s action. This momentum map is the composition
$$
\xymatrix{
T^*G \ar[r]^{\mm_H} \ar@/_9pt/[rr]_{\mm_S} & \h^* \ar@{->>}[r]^{ \iota_S^* } & {\mathfrak s}^*,
}
$$ where $\iota_S : {\mathfrak s} \to \h$ is the inclusion. Let us
identify $T^* G = G \times \g^*$ via left-translation. The momentum
map $\mm_{G_+}$ of $G$'s left action (resp. $\mm_{G_-}$, right
action) on $T^* G$ is equal, in this trivialization, to
$$\forall g\in G, x\in \g^*: \qquad \mm_{G_+}(g,x) = \Ad{g^{-1}}^* x, \qquad
\mm_{G_-}(g,x) = x.$$ 
In particular, $$\mm_{H} = \mm_{G_+}
\oplus - \mm_{G_-}.$$ 
The key property of the momentum map is its equivariance
$$\forall s\in S, P \in T^* G: \qquad \mm_S( s\cdot P ) = \Ad{s^{-1}}^*
\mm_S(P).$$

Let $\Sigma = G/U$, the quotient of $G$ by the free action of
$U$. The cotangent bundle of $\Sigma$ is known to be symplectomorphic
to the quotient $\mm_U^{-1}(0) / U$. Therefore, if $f_1,\ldots,f_n$ 
is a collection of $n=\dim \Sigma$ analytic functions on $T^*G$ that are
$U$-invariant, in involution, and functionally independent along
$\mm_S^{-1}(0)$, then they induce a real-analytically completely
integrable system on $T^* \Sigma$. In our case, we want $f_1$ to be a
convex hamiltonian; the natural choice is to have $f_1$ be the natural
bi-invariant metric on $T^* G$ induced by the Cartan-Killing form.

To construct the requisite integrals, it is necessary to digress
somewhat.

\subsubsection{Poisson algebras} \label{sssec:ddim}
Let $M$ be a smooth manifold. A Poisson structure on $M$ is a Lie
algebra structure $\pb{}{}$ on $\Ch{M}$ that makes $\Ch{M}$ into a Lie
algebra of derivations of $\Ch{M}$. Let $\A \subset \Ch{M}$ be a
subset--generally, a subalgebra--and denote by
\begin{equation} \label{eq:A}
\d\A_m = \spn{ \d_m f\ :\ f \in \A }, \hspace{10mm} Z(\A) = \left\{ f \in
\A\ :\ \pb{f}{g} \equiv 0 \ \  \forall g\in \A  \right\}.
\end{equation}
The subset $Z(\A)$ is the centre of $\A$: all functions in $\A$ are
integrals of the hamiltonian vector fields $X_f(\bullet) =
\pb{f}{\bullet} : f \in \A$. Let
\begin{equation} \label{eq:ddim}
\begin{split}
\ddim{\A}_m &= \dim \d\A_m,\\
\drank{\A}_m &= \dim \spn{ X_f(m)
  = \pb{f}{ \cdot }_m\ : \ f \in Z(\A) },
\end{split}
\end{equation}
be the differential dimension and the differential rank of $\A$ at
the point $m\in M$. We denote by $\ddim{\A}$ the maximum of $\left\{
\ddim{\A}_m\ :\ m \in M \right\}$ and by $\drank{\A}$ the maximum of
$\left\{ \drank{\A}_m\ :\  \ddim{\A}_m=\ddim{\A} \right\}$. If $\A$
contains a proper function, then $\drank{\A}+\ddim{\A}=\dim M$ implies
that the flow of any hamiltonian in the centre of $\A$ is integrable.

When $M = \h^*$ equipped with the canonical Poisson bracket, the
coadjoint orbits are symplectic leaves of the Poisson bracket. If a
subset $\A \subset \Cw{\h^*}$ satisfies
\begin{equation} \label{eq:inta}
\ddim{\A}_x = \frac{1}{2}\left( \dim H+ \rank H  \right),
\hspace{10mm} \drank{\A}_x = \frac{1}{2}\left( \dim H- \rank H  \right),
\end{equation}
for some regular $x\in \h^*$, then $\A| \orbit{y}$ defines a
real-analytically completely integrable system for all $y$ in an open
real-analytic subset of $\h^*$. If $\A$ satisfies the conditions
(\ref{eq:inta}), then $\A$ will be said to be completely integrable.

The centre of $\Cw{\h^*}$ is the set of real-analytic Casimirs; these
are functions that are constant on each coadjoint orbit. It is a
classic theorem of Cartan's that $Z(\Cw{\h^*})$ is finitely generated
by polynomials when $\h$ is semisimple. A technique that generates a
completely integrable algebra $\A$, that was discovered by
Mischenko and Fomenko~\cite{MF} and is related to Lax
representations, is the argument-shift technique.

\begin{qthm}[Mischenko-Fomenko 1978] \label{thm:MF}
Let $a \in \h^*$ and define $\A \subset \Cw{\h^*}$ by
\begin{equation} \label{eq:as}
\A := \left\{ f\ :\ \exists \lambda \in \R, g \in Z(\Cw{\h^*})\
\textrm{s.t.}\ f(x) = g(x+\lambda \cdot a)  \right\}.
\end{equation}
The algebra $\A$ is abelian and if $a$ is a regular element, then $\A$
is completely integrable.
\end{qthm}

Note that $\A$ always contains the Casimirs, and in particular, the
Cartan-Killing form of $\h^*$.

\subsubsection{The pull-back algebra } \label{sssec:pullback} 
Let $\A$ be a completely integrable subalgebra of $\Cw{\h^*}$ and let
$\B = \mm_H^* \A \subset \Cw{T^*G}$ be the algebra pulled-back to
$T^*G$. Since the momentum map $\mm_H$ is Poisson, $\B$ is also
abelian. We would like to compute $\drank{\B}$ and $\ddim{\B}$.

\begin{lemma} \label{lem:dind}
The differential rank and dimension of $\B$ equal $\dim
G$. Therefore, $\B$ defines a real-analytic completely integrable
subalgebra of $\Cw{T^*G}$.
\end{lemma}

\begin{proof}
Let $f \in \Cw{\h^*}$ and let $F = f \circ \mm_H$ be the pull-back of
$f$ by the momentum map. The chain rule shows that for $P \in T^* G$
\begin{equation} \label{eq:dPF}
\d_P F = \d_x f \cdot \d_P \mm_H,
\end{equation}
where $x=\mm_H(P)$. With the left-trivialization of $T^*G$, one can
write $P=(g,\mu)$ and a tangent vector $v\in T_P(T^*G)$ equals $\xi
\oplus \eta \in \g \oplus \g^*$. With this notation, and denoting
$\alpha_{\pm} = \pi_{\pm}(\d_x f)$ and $x_{\pm} = \jj_{\pm}^*(x)$, one
sees that $\d_PF \cdot v = 0$ iff
\begin{equation} \label{eq:dmm}
0 = \d_x f \cdot \d_P\mm_H \cdot v = \left\langle \alpha_+ ,
-\ad{\xi}^* x_+ + \Ad{g^{-1}}^* \eta \right\rangle - \left\langle
\alpha_- , \eta \right\rangle.
\end{equation}
If $\d_PF = 0$, then equation (\ref{eq:dmm}) vanishes for all
$v$. With $\eta=0$, this shows that $\alpha_+$ vanishes on
$\ad{\g}^*x_+$, which is the tangent space to $G$'s coadjoint orbit
through $x_+$. From this, equation (\ref{eq:dmm}) yields that 
\begin{equation} \label{eq:dfp}
\d_PF = 0 \qquad \iff \qquad \alpha_- = \Ad{g^{-1}}\alpha_+
\textrm{\ and\ } \left. \alpha_+ \right|_{\ad{\g}^*x_+}=0,
\end{equation}
and the second condition implies that $\left. \alpha_-
\right|_{\ad{\g}^*x_-}=0$, too. 

Assume that the point $x_+$ is a regular element in $\g^*_+$ (since
$x_-=-\mu$ and $x_+ = \Ad{g^{-1}}^*\mu$, one can equally assume that
$\mu$ is a regular point). In this case, the annihilator of
$\ad{\g}^*x_{\pm}$ is spanned by the derivatives of the Casimirs of
$\g_{\pm}^*$. Therefore, there is a Casimir $\phi$ of $\g_{+}^*$ such
that $\alpha_{+} = \d_{x_{+}}\phi$. The right-hand side of
(\ref{eq:dfp}) along with the fact that $x_-=-\Ad{g}^*x_+$ implies that
$\alpha_- = -\d_{x_-}\phi$. Therefore
\begin{equation} \label{eq:dfp1}
\d_PF = 0 \qquad \iff \qquad \exists \phi \in Z(\Cw{\g^*})\ \
\textrm{s.t.}\ \ \d_xf = \d_{x_+} \phi \oplus -\d_{x_-} \phi,  
\end{equation}
This implies that $\A$ contains the $1$-jets of all functions $F\in
\mm_H^*\Cw{\h^*}$ with $\d_PF=0$. Therefore, the dimension of $\d_P
\mm_H^*\A$ is equal to the dimension of $\d_x \A$ minus the dimension
of $\ker \d_P\mm_H^*$. The latter dimension equals $\rank \g$ by
(\ref{eq:dfp1}). Therefore,
\begin{equation} \label{eq:ddimA}
\ddim{\mm_H^* \A}_P = \frac{1}{2} \times ( \dim H + \rank{H}  ) -
\rank G = \dim G.
\end{equation}
Since $\B=\mm_H^*\A$ is abelian, this proves the lemma.
\end{proof}

\begin{thm} \label{thm:U}
Let $T \subset G$ be a maximal torus and assume that $U \subset T
\times T$ acts freely on $G$. Then, there are completely integrable
algebras $\A \subset \Cw{\h^*}$ such that $\B=\mm_H^*\A$ induces a
real-analytic, completely integrable convex hamiltonian on $T^*(G/U)$.
\end{thm}

\begin{proof}
$T^*(G/U)$ is canonically symplectomorphic to $\mm_U^{-1}(0)/U$, so it
  suffices to find an algebra $\A$ such that $\A$ is
  $\Ad{U}^*$-invariant and $\mm_U^{-1}(0)$ contains a regular point
  for $\B$.

To achieve $\Ad{U}^*$-invariance of $\A$, let the $a$ in the
Mischenko-Fomenko construction be chosen to lie in
$\t^* \oplus \t^*$. Since $T \times T$ is a maximal torus containing
$U$, there are regular elements $a \in \t^* \oplus \t^*$ and these are
stabilized by $\Ad{U}^*$. The equivariance of $\mm_H$ implies that
$\B$ is invariant under the action of $U$. This implies that $\B$ and
$\mm_U^*\Cw{\uu^*}$ Poisson commute; it also implies that $\B
| \mm_U^{-1}(0)$ induces a commutative Poisson algebra $\tilde\B$ of
real-analytic functions on $T^*(G/U)$.

The proof of Lemma~\ref{lem:dind} shows that $P=(g,\mu)\in T^*G$ is
regular for $\B$ if $\mu$ is a regular point in $\g^*$. Since $U$ is
contained in a maximal torus $T \times T$, the simple form of $\mm_H$
shows that the image of $\mm_H | \mm_U^{-1}(0) \cap T^*_1 G$ contains
the subspace of vectors $x=\eta \oplus -\eta$ such that $\eta \in
\t^{\perp}$. Since $\t^{\perp}$ contains regular elements, one
concludes that there are regular points for $\B$ in $\mm_U^{-1}(0)$.

Since $\B$ and $\mm_U^*\Cw{\uu^*}$ Poisson commute and each is
real-analytic, there is no loss in replacing $\B$ with the algebra $\B
+ \mm_U^*\Cw{\uu^*}$. To avoid a proliferation of notation, let $\B$
denote this expanded algebra. If $P \in \mm_U^{-1}(0)$ is a regular
point for $\B$, then we conclude that
\begin{equation} \label{eq:dindPU}
\begin{split}
\d_{U\cdot P} \tilde\B &\cong \d_P\B/\d_P \mm_U^*\Cw{\uu^*}\\
\implies \quad \ddim{\tilde\B}_{U\cdot P} &= \ddim{\B}-\dim U = \dim
G/U
\end{split}
\end{equation}
since the action of $U$ is free and $\B$ is completely
integrable. Therefore $\tilde\B$ is a real-analytic, completely
integrable algebra on $T^*(G/U)$. Since $\A$ contains the Casimirs of
$\h^*$, the Cartan-Killing form in $\A$ induces a riemannian metric
hamiltonian in $\tilde\B$. This is the convex hamiltonian that was sought.
\end{proof}

\subsection{Aloff-Wallach and Eschenburg Spaces } \label{ssec:esch} 

Recall the definition of an Eschenburg space~\cite{Esch}: let $U \cong
\U_1$ be a subgroup of $\SU_3 \times \SU_3$ such that the natural action of
$U$ on $\SU_3$ defined by
$$\forall u=(u_1,u_2) \in U, g \in \SU_3: \qquad u\cdot g = u_1 g
u_2^{-1}$$ is free. $U$ can be characterized in
terms of $4$ integers--$k,l,p,q$--as
$$U = \left\{  \diag(z^k,z^l,z^{-k-l}) \oplus \diag(z^p,z^q,z^{-p-q})\
:\ z \in S^1\right\}
$$
and $k,l,p,q$ satisfy
$$
\begin{array}{lllc}
\gcd(k-p,l-q), & \gcd(k-p,k+p+q), & \gcd(k+p+q,l-p),\\
\gcd(k-q,l-p), & \gcd(k-q,k+p+q), & \textrm{and } \gcd(k+p+q,l-q) & \textrm{equal } 1.\\
\end{array}
$$ Such $4$-tuples of integers are called admissible. Let $k,l,p,q$ be
an admissible quartet and let $U=U_{klpq}$ be such a group,
$M=M_{klpq}=\SU_3/U$. The manifold $M_{klpq}$ is an Eschenburg space.
When $k=l=0$, one has an Aloff-Wallach manifold~\cite{AW}.

Let $\kappa$ be a bi-invariant metric on $\SU_3$ and let $\kappa_M$ be
the submersion metric on $M$ induced by $\kappa$. Let $H$ and $H_M$ be
the induced fibre-quadratic hamiltonians on $T^*\SU_3$ and $T^*M$,
respectively.

\begin{thm} \label{thm:Esch}
The hamiltonian $H_M$ is real-analytically completely integrable for
any admissible quartet $k,l,p,q$.
\end{thm}

\begin{proof}
This is a simple corollary of the Theorem~\ref{thm:U}.
\end{proof}

\begin{remark} \label{rem:e}
 Here is a sketch of a pedestrian proof of theorem~\ref{thm:Esch}. The
functions $f_1,f_2,f_4,f_5$ defined in equation (\ref{al:f}) are
independent and in involution on $\su_3$; combined with $f_9= \det
\xi$, one obtains a completely integrable algebra of functions on
$\su_3$. The functions $f_1,f_2$ generate the coadjoint action of the
maximal torus of diagonal matrices. Therefore, the functions
$F_{i,\pm} = f_i \circ \mm_{G_{\pm}}$ with $G=\SU_3$ yield a
completely integrable algebra on $T^*G$ that is invariant under the
coadjoint action of the maximal torus in $G_+ \times G_-$ consisting
of diagonal matrices.

Paternain and Spatzier \cite{PS1} proved the integrability of $H_M$ on
Eschenburg spaces $M_{1,-1,2m,2m}$, using integrals like those in the
above paragraph along with some involved
computations. Bazaikin~\cite[Section 5]{Baz} proved the integrability
of a submersion geodesic flow on an Eschenburg space $M$ when
$M=H\backslash G/K$ where $G=\U_3 \oplus \U_2 \oplus \U_1$, $H=\U_2
\oplus\U_1$ and $K$ is isomorphic to $\U_1\oplus\U_1$ -- these
Eschenburg spaces are positively curved, in addition to having this
special bi-quotient structure. Theorem~\ref{thm:Esch} generalizes each
of these results.

Kruggel \cite{Kr} has obtained a homeomorphism and diffeomorphism
classification of Eschenburg spaces that satisfy his condition
$C$. Condition $C$ implies that the Eschenburg space is cobordant to a
union of 3 lens spaces. The second integral cohomology group of the
Eschenburg space $M=M_{klpq}$ is infinite cyclic with generator $u$;
the fourth integral cohomology group is finite cyclic of order $r$
with generator $u^2$. The first Pontryagin class $p_1(M)=p_1\cdot u^2$
and $u$ has a self-linking number $-s^{-1}/r \in \Q/\Z$, where
$s^{-1}$ is the inverse to $s$ in $\Z_r$. The final two invariants,
the Kreck-Stolz invariants $s_1, s_2 \in \Q/\Z$ are invariants of $M$
that Kruggel showed are calculable in terms of the eta-invariants of
the cobounding lens spaces. Each invariant is expressible in terms of
the parameters $k,l,p$ and $q$, but the formulae for $s_1, s_2$
involve transcendental functions. The invariants $r \in \Z$, $p_1,s
\in \Z_r$ and $s_2 \in \Q/\Z$ determine the oriented homeomorphism
type of $M_{klpq}$, while $r,p_1,s,s_1$ and $s_2$ determine the
oriented diffeomorphism type of $M_{klpq}$. Because $H^3(M;\Z_2)=0$,
there are a maximum of $28$ oriented smooth structures on a
topological Eschenburg space. In figure \ref{fig:eschspaces} one sees
the results of a numerical search for these $28$ oriented smooth
structures on the oriented topological Eschenburg space
$M_{-1,-1,-2,0}$. It is notable that, up to six decimal places, the
Kreck-Stolz invariant $s_1$ equals $i/28$ for $i=28,\ldots,1$.

\begin{center}
\begin{figure}[htb]
{
\setlength{\extrarowheight}{1pt}
\begin{tabular}{|rrrrl|rrrrl|}
\hline\hline
$k$ & $l$ & $p$ & $q$ & $s_1 \bmod 1$ & $k$ & $l$ & $p$ & $q$ & $s_1 \bmod 1$\\
\hline\hline
$-29$ &	$10$ &	$-28$ &	$6$ &	$1$       & 	$-21$ &	$-6$ &	$-18$ &	$-10$ &	$0.5$\\
$-38$ &	$-29$ &	$-66$ &	$22$ &	$0.964286$  &	$-5$ &	$-5$ &	$-6$ &	$-4$ &	$0.464286$\\
$-54$ &	$9$ &	$-52$ &	$4$ &	$0.928571$  &	$-13$ &	$2$ &	$-8$ &	$-6$ &	$0.428571$\\
$-17$ &	$-17$ &	$-18$ &	$-16$ &	$0.892857$  &	$-14$ &	$-5$ &	$-16$ &	$-2$ &	$0.392857$\\
$-6$ &	$-3$ &	$-8$ &	$0$ &	$0.857143$  &	$-9$ &	$-6$ &	$-12$ &	$-2$ &	$0.357143$\\
$-17$ &	$-14$ &	$-22$ &	$-8$ &	$0.821429$  &	$-38$ &	$-11$ &	$-48$ &	$8$ &	$0.321429$\\
$-14$ &	$-5$ &	$-18$ &	$2$ &	$0.785714$  &	$-22$ &	$-19$ &	$-40$ &	$12$ &	$0.285714$\\
$-1$ &	$-1$ &	$-2$ &	$0$ &	$0.75$    &	$-22$ &	$-1$ &	$-14$ &	$-12$ &	$0.25$\\
$-33$ &	$-6$ &	$-42$ &	$20$ &	$0.714286$  &	$-25$ &	$-1$ &	$-22$ &	$-6$ &	$0.214286$\\
$-46$ &	$-13$ &	$-32$ &	$-30$ &	$0.678571$  &	$-54$ &	$-9$ &	$-68$ &	$30$ &	$0.178571$\\
$-22$ &	$5$ &	$-20$ &	$0$ &	$0.642857$  &	$-39$ &	$-6$ &	$-32$ &	$-16$ &	$0.142857$\\
$-13$ &	$2$ &	$-14$ &	$6$ &	$0.607143$  &	$-29$ &	$-14$ &	$-32$ &	$-10$ &	$0.107143$\\
$-38$ &	$-11$ &	$-40$ &	$-8$ &	$0.571429$  &	$-11$ &	$1$ &	$-12$ &	$4$ &	$0.071429$\\
$-22$ &	$-1$ &	$-26$ &	$12$ &	$0.535714$  &	$-9$ &	$-9$ &	$-10$ &	$-8$ &	$0.035714$\\
\hline
\end{tabular}
}
\caption{Representative Eschenburg spaces $M_{klpq}$ in the oriented
  homeomorphism class with $r=1$, $s=0$, $p_1=0$ and $s_2=0.25 \bmod
  1$. The Kreck-Stolz invariant $s_1$ has been rounded to six decimal
  places. See \cite{Butler:2008a} for details of the
  computations.} \label{fig:eschspaces}
\end{figure}
\end{center}

\end{remark}

\section{ Conclusion }
In~\cite{DRWCT}, Dullin, Robbins, Waalkens, Creagh and Tanner
demonstrate that a cohomologically non-vanishing Maslov cocycle
constrains the monodromy of a completely integrable
system. Specifically, they show that the cohomology class of the
Maslov cocycle, if non-zero on a lagrangian torus, is a common
eigenvector of the monodromy group of the lagrangian fibration. Their
work {\em assumes} the non-triviality of the Maslov cocycle, and ends
with the question: {\em Does the Maslov cocycle of an invariant torus
of a natural mechanical hamiltonian on $T^*\R^n$ always vanish?}

The answer to their question is {\em yes} and was proven by Viterbo in
the work cited above \cite{Viterbo}. Viterbo proved that if $F
\subset T^*\R^n$ is a lagrangian torus, then there is a cycle
on $F$ whose Maslov index is an even integer between $2$ and $n+1$
inclusive. The proof uses Conley-Zehnder theory, and the strength of
the result is the constraint on how far the Maslov cocycle may be from
primitive. This is used to prove the above-cited result on the
vanishing of the Maslov cocycle when $\rho$ has a non-zero degree.

There are several natural questions that arise from the note of Dullin,
{\it et. al.}. First, there are higher-dimensional Maslov cocycles
that measure the higher singularities of $\rho$; the cohomological
non-triviality of these cocycles further constrains the monodromy of a
completely integrable system. What is it possible to say about their
non-vanishing? Second, if $F \subset T^* \Sigma$ is a lagrangian torus
whose Maslov class $\iota_F^*(d)$ is cohomologically non-trivial, must
this Maslov class be close to primitive? That is, how far does
Viterbo's results generalize? We note that Viterbo himself has
obtained one generalization~\cite{Viterbo2} and that Bialy
\cite{Bialy} has shown in two degrees of freedom that the Maslov class
is twice a primitive element in many cases.


\def\cprime{$'$} \def\cprime{$'$}

\end{document}